\documentclass[10pt]{article}
\usepackage[utf8]{inputenc}

\usepackage{amsmath,latexsym}
\usepackage[psamsfonts]{amssymb}
\usepackage[bookmarks,bookmarksnumbered,pdfstartview=FitBH,backref=page]{hyperref}
\hypersetup{
colorlinks=true, 
citecolor= blue,
linktoc= all, 
linkcolor= red, 
urlcolor= green,
}
\usepackage[backref=page]{hyperref}

\usepackage{theorem}
\theoremstyle{changebreak}
\newtheorem{theorem}{Theorem}[section]
\newtheorem{proposition}[theorem]{Proposition}
\newtheorem{question}[theorem]{Question}
\newtheorem{corollary}[theorem]{Corollary}
\newtheorem{claim}[theorem]{Claim}

\def\defin#1{\textbf{\emph{#1}}}
\newcommand\alert[1]{\textbf{#1}}
\def\CQFD{\nobreak\hfill $\blacksquare$ \goodbreak}

\newcommand\Nmath{\mathbb{N}}
\newcommand\Rmath{\mathbb{R}}
\newcommand\Zmath{\mathbb{Z}}
\newcommand\Qmath{\mathbb{Q}}
\newcommand\Smath{\mathbb{S}}
\newcommand\Hmath{\mathbb{H}}

\newcommand\FF{{\bf F}}

\newcommand\Aut{\mathrm{Aut}}
\newcommand\Out{\mathrm{Out}}
\newcommand\MCG{\mathrm{MCG}}
\newcommand\GL{\mathrm{GL}}
\newcommand\SL{\mathrm{SL}}
\def\PSL{\mathrm{PSL}}

\DeclareMathOperator{\acting}{\curvearrowright}
\DeclareMathOperator\vol{\mathrm{Vol}}
\DeclareMathOperator\image{im}

\newcommand\RR{{\mathcal{R}}}
\renewcommand\SS{{\mathcal{S}}}

\begin{document}

\title{On the top-dimensional $\ell^2$-Betti numbers}
\author{Damien Gaboriau\footnote{CNRS \& U.M.P.A., Ecole Normale Supérieure de Lyon, UMR 5669, 46 all{\'e}e d'Italie, 69364 Lyon cedex 07, France.},
Camille Noûs\footnote{Laboratoire Cogitamus; https://www.cogitamus.fr/}}

\maketitle

\begin{abstract}{The purpose of this note is to introduce a trick which relates the (non)-vanishing of the top-dimensional $\ell^2$-Betti numbers of actions with that of sub-actions.
We provide three different types of applications: we prove that the $\ell^2$-Betti numbers of $\Aut(\FF_n)$ and $\Out(\FF_n)$ (and of their Torelli subgroups) do not vanish in degree equal to their virtual cohomological dimension, we prove that the subgroups of the $3$-manifold groups have vanishing $\ell^2$-Betti numbers in degree $3$ and $2$ and we figure out the ergodic dimension of certain direct products of the form $H\times A$ where $A$ is infinite amenable.
}\end{abstract}

\noindent
\begin{center}
\textbf{Résumé en Français} 
\medskip

\begin{minipage}{0.85\linewidth}
{\small Le but de cette note est d'introduire une astuce qui relie l'annulation (ou la non-annulation) du nombre de Betti $\ell^2$ en dimension maximale des actions d'un groupe avec l'annulation pour ses sous-actions. On fournit trois différents types d'applications :
on montre que les nombres de Betti $\ell^2$ de $\Aut(\FF_n)$ et $\Out(\FF_n)$ (et de leurs sous-groupes de Torelli) ne s'annulent pas en degré égal à leur dimension cohomologique virtuelle ; on prouve qu'un sous-groupe quelconque du groupe fondamental d'une variété compacte de dimension $3$ a ses nombres de Betti $\ell^2$ nuls en degré $3$ et $2$
et enfin, on parvient à déterminer la dimension ergodique de certains produits directs de la forme $H\times A$ où $A$ est moyennable infini.}
\end{minipage}
\end{center}

\medskip

\noindent
\textbf{keywords:}
{$\ell^2$-Betti numbers, measured group theory, cohomological dimension, ergodic dimension, $\Out(\FF_n), \Aut(\FF_n)$, $3$-dimensional manifolds.}

\noindent
\textbf{MSC:}
{37A20, 19K56, 20F28, 20E15, 57Mxx.}

\vskip 20pt

\section{Presentation of the results}
The $\ell^2$-Betti numbers were introduced by Atiyah \cite{Ati76}, in terms of heat kernel, for free cocompact group actions on manifolds and were extended to the framework of measured foliations by Connes \cite{Connes-79}. They acquired the status of group invariants thanks to Cheeger and Gromov \cite{CG86} who provided us with the definition of the $\ell^2$-Betti numbers of an arbitrary countable group~$\Gamma$:
\begin{equation*}
\beta_0^{(2)}(\Gamma), \beta_1^{(2)}(\Gamma), \beta_2^{(2)}(\Gamma), \cdots, \beta_k^{(2)}(\Gamma), \cdots
\end{equation*} 
Their extension to standard probability measure preserving actions and equivalence relations by the first author \cite{Gab02} opened the connection with the domain of orbit equivalence, offering in return some general by-products, for instance the \alert{$\ell^2$-proportionality principle} \cite[Corollaire 0.2]{Gab02}:
{\em If $\Gamma$ and $\Lambda$ are lattices in a locally compact second countable (lcsc) group $G$ with Haar measure $\vol$, then their $\ell^2$-Betti numbers are related as their covolumes:
$\frac{\beta_k^{(2)}(\Gamma)}{\vol (\Gamma\backslash G)}=\frac{\beta_k^{(2)}(\Lambda)}{\vol (\Lambda\backslash G)}$}.

Over the years, the $\ell^2$-Betti numbers have been proved to provide very useful invariants in geometry, in $3$-dimensional manifolds,  in ergodic theory,  in operator algebras and in many aspects of discrete group theory such as  geometric, resp. measured, resp. asymptotic group theory. 
We refer to \cite{Eck00} for an introduction to the subject and to the monographies \cite{Luc02,Kammeyer-book-intro-L2}.

\medskip
The term top-dimension used in the title may have different meanings.
At first glance, we mean the dimension of some contractible simplicial complex on which our group $\Gamma$ acts simplicially and properly (i.e., with finite stabilizers).
For the purpose of computing $\ell^2$-Betti numbers, one can consider the action of some finite index subgroup of $\Gamma$. In many interesting cases, the group $\Gamma$ is indeed virtually torsion-free. Then, the \defin{virtual geometric dimension} (the minimal dimension of a contractible simplicial  complex on which a finite index subgroup acts simplicially and freely) can be used as a better (i.e., lower) top-dimension for $\Gamma$. Observe that the $\ell^2$-Betti numbers must vanish in all degrees  above this dimension. In view of the Eilenberg-Ganea Theorem \cite{Eilenberg-Ganea-1957} (see also \cite[Chapter VIII.7]{Bro83}), if the \defin{virtual cohomological dimension (vcd)} of $\Gamma$ is finite and greater than three then it coincides with the virtual geometric dimension. 
The vanishing or non-vanishing of $\ell^2$-Betti numbers in some degree is an invariant for lattices in the same lcsc group (as the $\ell^2$-proportionality principle above indicates), and it is more generally an invariant 
of measure-equivalence \cite[Théorème 6.3]{Gab02}. In contrast, the  virtual cohomological dimension is not: for instance cocompact versus non-cocompact lattices in $\SL(d,\Rmath)$ have different vcd. This nominates the \defin{ergodic dimension} as a better notion of top-dimension.
This is intrinsically an invariant of measured group theory introduced in \cite[Définition 6.4]{Gab02} (see Section~\ref{sect: proofs with erg dim} and also \cite{Gab-erg-dim}) which mixes geometry and ergodic theory. It is bounded above by the virtual geometric dimension and is often much less. Our trick (Theorems~\ref{th: top dim L2 pmp equiv rel, geom dim} and \ref{th: top dim L2 R-complex}) also applies to it.

\subsection{\texorpdfstring{$\Aut(\FF_n)$ and $\Out(\FF_n)$}{Aut Fn and Out Fn}}
While the $\ell^2$-Betti numbers of many classic groups are quite well understood, this is far from true for the groups $\Aut(\FF_n)$ and $\Out(\FF_n)$
of automorphisms (resp. outer automorphisms) of the free group $\FF_n$ on $n\geq 3$ generators.
 These groups share many algebraic features with both the group $\GL(n,\Zmath)$ and with the mapping class group $\MCG(S_{g})$ of the surface $S_g$ of genus $g$. 
 One reason is that all these groups are  (outer) automorphism groups of the most primitive discrete groups ($\FF_n$, $\Zmath^n$ and $\pi_1(S_{g})$ respectively) and the three families begin with the same group $\Out(\FF_2)\simeq\GL(2,\Zmath)\simeq  \MCG(S_{1})$.
 These empirical similarities have served as guiding lines for their study, see for instance \cite{CV86,Bridson-Vogtmann-2006,Vogtmann-2006-ICM-survey}.

By the work of Borel \cite{Bor85}, the $\ell^2$-Betti numbers of the cocompact lattices of $\GL(n,\Rmath)$ are known to all vanish when $n\geq 3$. 
The same holds for the non-cocompact ones like $\GL(n,\Zmath)$ by the $\ell^2$-proportionality principle. 
The mapping class group $\MCG(S_{g})$ is virtually torsion-free, 
  and when $g>1$, all its $\ell^2$-Betti numbers vanish except in degree equal to the middle dimension $3g-3$ of its Teichm\"uller space (see for instance \cite[Appendix D]{Kida-08}). These behaviors are very common for $\ell^2$-Betti numbers of the classic groups: most of them vanish, and when a non-vanishing  happens it is only in the middle dimension of ``the associated symmetric space''.

Culler-Vogtmann \cite{CV86} invented the Outer space $CV_n$ as an analogue of the Teichm\"uller space in order to transfer (rarely straightforwardly) 
 the geometric techniques of  Thurston for the mapping class groups to $\Out(\FF_n)$. It is also often thought of as an analogue of the symmetric space of lattices in Lie groups.
It has dimension $3n-4$ and admits an $\Out(\FF_n)$-equivariant deformation retraction onto a proper contractible simplicial complex, the spine of the outer-space, of dimension $2n-3$ which is thus exactly the virtual cohomological dimension of $\Out(\FF_n)$ \cite[Corollary 6.1.3]{CV86} (a lower bound being easy to obtain). An avatar of $CV_n$ can be used to show that the virtual cohomological dimension of $\Aut(\FF_n)$ is  $2n-2$ \cite[pp. 59-61]{Hatcher-1995-Homol-stab-aut-free-gps}.

\begin{theorem}\label{th: beta 2n-3 Out Fn}
The $\ell^2$-Betti numbers of the groups $\Out(\FF_n)$ and $\Aut(\FF_n)$ ($n\geq 2$) do not vanish in degree equal to their virtual cohomological dimensions $2n-3$ (resp. $2n-2$): 
\begin{equation*}
\beta_{2n-3}^{(2)}(\Out(\FF_n))>0 \text{ and } \beta_{2n-2}^{(2)}(\Aut(\FF_n))>0.
\end{equation*}
\end{theorem}
The rational homology of $\Out(\FF_n)$ is  very intriguing. It was computed explicitly using computers by Ohashi \cite{Ohashi-2008-rat-homol-Out-Fn} up to $n=6$.
Then Bartholdi \cite{Bartholdi-2016-homol-out-F7} proved for $n=7$ that $H_k(\Out(\FF_7);\Qmath)$ is trivial except for $k=0,8,11$, when it is  $1$-dimensional.
The non-zero classes for $k=8,11$ were a total surprise, since they are not generated by Morita classes.
Moreover, the rational homology of both $\GL(n,\Zmath)$ and $\MCG(S_g)$ vanishes in the virtual cohomological dimension, and everyone expected the same would be true for $\Out(\FF_n)$.   
In view of the L{\"u}ck approximation \cite{Luck-1994-approx},  Theorem~\ref{th: beta 2n-3 Out Fn} implies that in degree equal to their vcd, 
the rational homology grows indeed linearly along towers. More precisely, these groups being residually finite \cite{Baumslag-1963-Autom-gps-resid-finite,Grossman-74}, for every sequence of finite index normal subgroups $(\Gamma_i)_i$ which is decreasing with trivial intersection in $\Out(\FF_n)$ (resp. $\Aut(\FF_n)$), then
\begin{equation*}
\lim\limits_{i\to \infty} \frac{\dim H_{2n-3}(\Gamma_i;\Qmath)}{[\Out(\FF_n):\Gamma_i]}>0, \ \ \text{ resp. } \ \ \ \lim\limits_{i\to \infty} \frac{\dim H_{2n-2}(\Gamma_i;\Qmath)}{[\Aut(\FF_n):\Gamma_i]}>0.
\end{equation*}
The mystery top-dimensional classes implicitly exhibited here for large finite index subgroups ``come'' from a poly-free subgroup $\FF_2\ltimes \FF_2^{2n-4}$ of $\Out(\FF_n)$. 
In a work in progress with Laurent Bartholdi, we build on this remark to produce more explicit classes \cite{Bar-Gab}.
We also work on discovering other $\ell^2$-Betti numbers for $\Out(\FF_n)$.
Results of Smillie and Vogtmann suggest that the (rational) Euler characteristic (equivalently the standard Euler characteristic of any torsion-free finite index subgroup) of $\Out(\FF_n)$ should always be negative and this has been indeed proved very recently by Borinsky and Vogtmann \cite{Borinsky-Vogtmann-2019}.
A positive answer to the following question would deliver another demonstration.
\begin{question}
Do all $\beta_{k}^{(2)}(\Out(\FF_n))$ for $k\not=2n-3$ vanish?
\end{question}
Theorem~\ref{th: beta 2n-3 Out Fn}
will be proved in section~\ref{sect: proof th Out Fn Aut Fn}.

The canonical homomorphisms of $\Aut(\FF_n)$ and $\Out(\FF_n)$ to $\GL(n,\Zmath)$ lead to the short exact sequences
\begin{eqnarray}
1\to {\mathcal{T}}_n \to \Out(\FF_n)\overset{\phi_n}{\to} \GL(n,\Zmath)\to 1. \label{eq Nielsen-Magnus morph Out}\\
1\to {\mathcal{K}}_n \to \Aut(\FF_n)\overset{}{\to} \GL(n,\Zmath)\to 1. \label{eq Nielsen-Magnus morph Aut}
\end{eqnarray}
The left hand side groups ${\mathcal{T}}_n$ and ${\mathcal{K}}_n$, called the \defin{Torelli groups}, have cohomological dimension $2n-4$ and $2n-3$ \cite{Bestvina-Bux-Margalit-2007}.

\begin{theorem}\label{th: top dim beta 2 of Torelli for Out Fn and Aut Fn}
The $\ell^2$-Betti numbers of the Torelli groups ${\mathcal{T}}_n$ and ${\mathcal{K}}_n$ ($n\geq 2$) do not vanish in degree equal to their virtual cohomological dimensions $2n-4$ (resp. $2n-3$): 
\begin{equation*}
\beta_{2n-4}^{(2)}({\mathcal{T}}_n)\not= 0 \text{ and } \beta_{2n-3}^{(2)}({\mathcal{K}}_n)\not= 0.
\end{equation*}
\end{theorem}
This is proved in Section~\ref{sect: Torelli subgroup}.

\subsection{Fundamental groups of compact manifolds}

We now switch to another type of application.
This one necessitates the full strength of the measured framework of Theorem~\ref{th: top dim L2 pmp equiv rel, geom dim} below.
The (virtual) cohomological dimension of the fundamental group $\pi_1(M)$ of a compact aspherical $d$-dimensional manifold $M$ is clearly $\leq d$, with equality when $M$ is closed. However, with Conley, Marks and Tucker-Drob we sharpened this in \cite{Conley-Gaboriau-Marks-Tucker-Drob-2016} by showing that $\Gamma=\pi_1(M)$ has ergodic dimension $\leq d-1$. 
This means that with the help of an auxiliary probability measure preserving  free $\Gamma$-action, one gains one on the top-dimension (see Section \ref{sect: proof th manifolds}).
And of course the smaller the ergodic dimension, the better the top-dimension. Thus the importance of Questions~\ref{quest: erg dim asph mfld groups}.
So far, we obtain:
\begin{theorem}\label{th: d-mfld and top L2 Betti}
Let $\Gamma$ be the fundamental group of a compact connected \alert{aspherical} manifold $M$ of dimension $d\geq 3$. Let $\Lambda\leq\Gamma$ be any subgroup.
Then $\beta_{d}^{(2)}(\Lambda)=0$.
If moreover $\beta_{d-1}^{(2)}(\Gamma)=0$ then $\beta_{d-1}^{(2)}(\Lambda)=0$.
 \end{theorem}
Of course all the  $\ell^2$-Betti numbers of $\Lambda$ vanish in degree $>d$.
Observe that the asphericity is a necessary condition in this statement since for instance $\FF_2^{4}$ is the fundamental group of some compact $4$-manifold while its $4$-th $\ell^2$-Betti number equals $1$.
Recall that the Singer Conjecture predicts that {\em the $\ell^2$-Betti numbers of a closed aspherical manifold $M$  are concentrated in the middle dimension,
i.e., if $\beta_{k}^{(2)}(\pi_1(M))>0$ then $2k=$ the dimension of $M$.} The "moreover" assumption of Theorem~\ref{th: d-mfld and top L2 Betti} would then be satisfied automatically.
The Singer Conjecture holds in particular for closed hyperbolic manifolds \cite{Dodziuk-1979-harmon-sym-mflds}.
Given the recent progress on $3$-dimensional manifolds (\cite{Perelman-2002-entropy-Ricci,Perelman-2003-Ricci-surgery}, see also \cite{Kleiner-Lott-2008-Notes-Perelman,BBBMP-2010}), we obtain a more general statement:
\begin{theorem}\label{th: 3-dim manifold and L2 Betti}
Let $\Gamma$ be  the  fundamental group of a 
connected compact $3$-dimensional manifold.
The $\ell^2$-Betti numbers of any subgroup $\Lambda\leq\Gamma$ vanish in all degrees $k\geq 2$:
\begin{equation*}
\beta_k^{(2)}(\Lambda)=0.
\end{equation*}
In particular, if $\Lambda$ is infinite then $\chi^{(2)}(\Lambda)\in [-\infty, 0]$.
\end{theorem}
Here $\chi^{(2)}(\Lambda)$ is the $\ell^2$-Euler characteristic of $\Lambda$. It coincides with the virtual Euler characteristic when the latter is defined.
Observe that the $3$-manifold in this theorem can have boundary, can be non-orientable and is not necessarily aspherical.
While the vanishing in degree $3$ for subgroups could have been expected, it is more surprising in degree $2$.
These results are proved in Section~\ref{sect: proof th manifolds}.

\subsection{Ergodic dimension}
Let's now switch to the third type of applications.
The non-vanishing of the $\ell^2$-Betti number in some degree $d$ for some subgroup $\Lambda$ of a countable group $\Gamma$ promotes clearly $d$ to a lower bound of the virtual geometric dimension of $\Gamma$. Although the ergodic dimension is bounded above by the virtual geometric dimension, 
$d$ is even a lower bound of the ergodic dimension of $\Gamma$
\cite[Corollaire 3.17, Corollaire 5.9]{Gab02}.
In case $\beta_d^{(2)}(\Gamma)=0$, then $d+1$ is upgraded a lower bound:
\begin{theorem}
\label{th: L2 erg-dim groups}
If $\Gamma$ is a countable discrete group of ergodic dimension (resp. virtual geometric dimension) $\leq d$ and if $\Lambda\leq\Gamma$ is any subgroup such that $\beta_d^{(2)}(\Lambda)\not=0$, then $\beta_d^{(2)}(\Gamma)\not=0$ and the ergodic dimension of $\Gamma$ is $d$.
\end{theorem}
This statement is an immediate application of Theorem~\ref{th: top dim L2 R-complex}.
It is worth recalling a result in this spirit: {\em  If $\Gamma$ is non-amenable and satisfies $\beta_1^{(2)}(\Gamma)=0$ then its ergodic dimension is $\geq 2$ \cite[Prop. 6.10]{Gab02}}. The non-amenability assumption plays here the role of a subgroup with non-zero $\beta_1^{(2)}$. And this is not just an analogy since non-amenable groups contain, in a measurable sense, a free subgroup $\FF_2$ \cite{Gaboriau-Lyons-2009}.

As a corollary, one computes the ergodic dimension of such groups as 
 $\FF_2^{d}\times \Zmath$: it is $d+1$.
 As another example $\Out(\FF_n)\times \Zmath^{k}$ (resp.  $\Aut(\FF_n)\times \Zmath^{k}$)  has ergodic dimension $2n-2$ (resp. $2n-1$).
 More generally, 
 \begin{corollary}\label{cor: ergim and L2}
 If $\Lambda$ has ergodic dimension $d$ and $\beta_d^{(2)}(\Lambda)\not=0$, then for any infinite amenable group $B$, the direct sum $\Lambda\times  B$ has ergodic dimension $d+1$.
 \end{corollary}
All the $\ell^2$-Betti numbers of $\Lambda\times B$ equal $0$. Observe that the condition $\beta_d^{(2)}(\Lambda)\not=0$ is necessary since for instance 
 $(\Lambda\times B)\times B=\Lambda\times  (B\times  B)$ has also ergodic dimension $d+1$.

\subsection{A top-dimensional $L^2$-Betti number result}
The different statements announced above use at some point variants of the general trick (Theorem~\ref{th: top dim L2 R-complex}) involving a probability measure preserving standard equivalence relation $\RR$ with countable classes (\defin{pmp equivalence relation} for short), a standard sub-relation $\SS$ and a simplicial discrete $\RR$-complex together with their $L^2$-Betti numbers\footnote{Observe the debatable use introduced in \cite{Gab02} of capital letter $L^2$ for equivalence relations versus the cursive lowercase $\ell^2$ for groups.}; see sections~\ref{sect: proof of R-complex top-dim trick} and \ref{sect: proofs with erg dim} where the notions are recalled. %
The specialization of Theorem~\ref{th: top dim L2 R-complex} to proper actions (simplicial actions with finite stabilizers) which is appropriate for geometric dimension will be given its own proof  in section~\ref{sect: top-dim trick groups} for the reader's convenience and as a warm-up to section~\ref{sect: proof of R-complex top-dim trick}.
Let's denote by $\beta_d^{(2)}(\Gamma\acting L)$ the $d$-th $\ell^2$-Betti number of the action of $\Gamma$ on $L$, also denoted by countless different manners in the literature 
such as $\beta_d(L,\Gamma)$, $\beta_d^{(2)}(L,\Gamma)$, $\beta_d^{(2)}(L:\Gamma)$ or $b^d_{(2)}(L:\Gamma)$.
\begin{theorem}[Proper actions version]
\label{th: top-dim groups}
Let $\Gamma$ be a countable discrete group and $\Lambda\leq\Gamma$ be a subgroup.
If $\Gamma\acting L$ is a proper action on a $d$-dimensional simplicial complex such that the restriction to $\Lambda$ satisfies $\beta_d^{(2)}(\Lambda\acting L)\not=0$, then $\beta_d^{(2)}(\Gamma\acting L)\not=0$.
\end{theorem}

\medskip
Specializing Theorem~\ref{th: top dim L2 R-complex} to a \alert{contractible} $\RR$-complex, one obtains a statement involving the $L^2$-Betti numbers of the pmp equivalence relation \cite[Théorème 3.13, Définition 3.14]{Gab02} and of its sub-relations. The minimal dimension of such a contractible complex defines the \defin{geometric dimension} of $\RR$ (see the proof of Theorem~\ref{th: top dim L2 pmp equiv rel, geom dim}).
\begin{theorem}[Geometric dimension of pmp equivalence relation]
\label{th: top dim L2 pmp equiv rel, geom dim}
If $\RR$ is a pmp equivalence relation on the standard space $(X,\mu)$  of geometric dimension $\leq d$ for which the $L^2$-Betti number in degree $d$ vanishes ($\beta_d^{(2)}(\RR,\mu)=0$) then every standard sub-equivalence relation $\SS\leq\RR$ satisfies $\beta_d^{(2)}(\SS,\mu)=0$.
\end{theorem}

\subsection*{Acknowledgments} The authors would like to thank Michel Boileau and Jean-Pierre Otal for enlightening discussions on $3$-dimensional manifolds. They also thank Clara L{\"o}h and Karen Vogtmann for their valuable comments on a preliminary version. 
Camille Noûs embodies the collective's constitutive role in the creation of scientific knowledge, which is always part of a collegial organization and the legacy  of previous work.

This work is supported by the ANR project GAMME (ANR-14-CE25-0004), by the CNRS and by the LABEX MILYON (ANR-10-LABX-0070) of Université de Lyon, within the program ``Investissements d'Avenir" (ANR-11-IDEX-0007) operated by the French National Research Agency (ANR).

\section{Proof of Theorem~\ref{th: top-dim groups} on proper simplicial actions}
\label{sect: top-dim trick groups}

Recall \cite[(2.8) p. 198]{CG86} (see also \cite[Section 1.2]{Gab02}) that for a proper non-cocompact action $\Gamma\acting L$, the $\ell^2$-Betti numbers are defined as follows:
Consider any increasing exhausting sequence $(L_i)_{i\in \Nmath}$ of cocompact $\Gamma$-invariant subcomplexes of $L$. For each dimension $k$, for each $i\leq j$, the inclusion $L_i\subset L_j$ induces a $\Gamma$-equivariant map $\bar{H}_k^{(2)}(L_i)\to \bar{H}_k^{(2)}(L_j)$ between the reduced $\ell^2$-homology spaces. The von Neumann $\Gamma$-dimension of the closure of the image $\image \left(\bar{H}_k^{(2)}(L_i)\to \bar{H}_k^{(2)}(L_j)\right)$ is decreasing in $j$ and increasing in $i$. The $k$-th $\ell^2$-Betti number of the action is defined as:
\begin{equation*}
\beta_k^{(2)}(\Gamma\acting L):=
\lim\limits_{i\to \infty} \nearrow \lim\limits_{
\begin{smallmatrix} j\to \infty\\ j\geq i
\end{smallmatrix}}\searrow
\, \dim_{\Gamma}\left(\overline{\image} \left(\bar{H}_k^{(2)}(L_i)\to \bar{H}_k^{(2)}(L_j)\right)\right).
\end{equation*}
This is easily seen to be independent of the choice of the exhausting sequence.
The  $k$-th $\ell^2$-Betti number of the group $\Gamma$ is defined as the $k$-th $\ell^2$-Betti number $\beta_k^{(2)}(\Gamma\acting L)$ for any proper contractible (or even only $k$-contractible) $\Gamma$-complex $L$ and this is independent of the choice of $L$.

The key observation is that for any  $d$-dimensional complex $M$ the reduced $\ell^2$-homology, defined from the $\ell^2$-chain complex 
\begin{equation*}
0\longleftarrow C_{0}^{(2)}(M)\overset{\partial_1}{\longleftarrow}  C_{1}^{(2)}(M) \cdots \longleftarrow C_{d-1}^{(2)}(M)\overset{\partial_d}{\longleftarrow}  C_{d}^{(2)}(M)\longleftarrow 0
\end{equation*}
boils down in dimension $d$ to the kernel of the boundary map 
\[\bar{H}_{d}^{(2)}(M)=H_{d}^{(2)}(M)=\ker \partial_d^{M}:=\ker \left(C_d^{(2)}(M)\overset{\partial_{d}}{\to} C_{d-1}^{(2)}(M)\right).\]
Of course, for the boundary operators to be bounded, 
$M$ needs here to have \defin{bounded geometry}, i.e., it admits a uniform bound on the \defin{valencies} (the number of simplices a vertex belongs to).

Since the  injective maps induced on $\ell^2$-chains by the inclusions $L_i\subset L_j$ commute with boundaries, it follows that 
\begin{eqnarray}
\beta_d^{(2)}(\Gamma\acting L)&:=&\lim\limits_{i\to \infty}\lim\limits_{
\begin{smallmatrix} j\to \infty\\ j\geq i
\end{smallmatrix}
}\, \dim_{\Gamma}\image \left(\ker \partial_{d}^{L_i} \lhook\joinrel\longrightarrow  \ker \partial_{d}^{L_j}\right)
\\
&=&\lim\limits_{i\to \infty}\nearrow\dim_{\Gamma}\ker \partial_{d}^{L_i}.\label{eq:beta top as kernel}
\end{eqnarray}

Consider, for the restricted action $\Lambda\acting L$, an increasing exhausting sequence $(K_i)_{i\in \Nmath}$ of cocompact $\Lambda$-invariant subcomplexes of $L$. By assumption, for $i$ large enough, $\dim_{\Lambda}\ker \partial_{d}^{K_i}\not=0$, so that $\ker \partial_{d}^{K_i}\not=\{0\}$. Let $L_i:=\cup_{\gamma\in \Gamma} \gamma K_i$ be the $\Gamma$-saturation of the $K_i$. It is $\Gamma$-invariant and $\Gamma$-cocompact. Again by commutation with boundaries of the injective maps induced on $\ell^2$-chains by the inclusion $K_i\subset L_i$, we also have $\ker \partial_{d}^{L_i}\not=\{0\}$. The $\Gamma$-saturations $L_j$ of the $K_j$ give an
increasing exhausting sequence $(L_j)_{j\in \Nmath}$ of cocompact $\Gamma$-invariant subcomplexes of $L$.
In view of formula (\ref{eq:beta top as kernel}) and since the von Neumann dimension is faithful, we have
$\beta_d^{(2)}(\Gamma\acting L)\not=0$ .\CQFD

\section{\texorpdfstring{Proof of Theorem~\ref{th: beta 2n-3 Out Fn} that $\beta_{2n-3}^{(2)}(\Out(\FF_n))>0$}{Proof th that beta 2n-2 Out Fn non zero}}
\label{sect: proof th Out Fn Aut Fn}

We begin by recalling what is known about the $\ell^2$-Betti numbers of $\Aut(\FF_n)$ and $\Out(\FF_n)$.
The groups fit into a canonical short exact sequence 
\begin{eqnarray}
1\to \FF_n \to \Aut(\FF_n) \overset{\theta_n}{\to}  \Out(\FF_n)\to 1. \label{eq: Aut to Out}
\end{eqnarray}
When $n=2$, the group $\Out(\FF_2)\simeq \GL(2,\Zmath)$ admits a single non-vanishing $\ell^2$-Betti number, namely $\beta_1^{(2)}(\GL(2,\Zmath))=1/24$ in degree $1$, exactly the middle dimension of its associated symmetric space and also the virtual geometric dimension of $\GL(2,\Zmath)$.
It follows that $\Aut(\FF_2)$ has an index $24$ subgroup isomorphic with $\FF_2\ltimes \FF_2$ so that its $\ell^2$-Betti numbers vanish except $\beta_2^{(2)}(\Aut(\FF_2))=1/24$ (see for instance Proposition~\ref{prop: L2 Betti poly-free groups}).
When $n\geq 3$, the kernel ${\mathcal{T}}_n$ of $\phi_n$ (sequence \eqref{eq Nielsen-Magnus morph Out}) is a finitely generated infinite normal subgroup of infinite index by \cite{Nielsen-1924,Magnus-1935} (and clearly the same holds for the kernel $\FF_n$ of $\theta_n$ (sequence \eqref{eq: Aut to Out}). It follows that $\beta_1^{(2)}(\Out(\FF_n))=\beta_1^{(2)}(\Aut(\FF_n))=0$ by using for instance \cite[Theorem  3.3 (5)]{Luc98b}  or  \cite[Théorème 6.8]{Gab02}:
{\em The middle group $H$ of a short exact sequence $1\to N\to H\to Q\to 1$ of infinite groups has $\beta_1^{(2)}(H)=0$ as soon as $\beta_1^{(2)}(N)<\infty$}. 
Remark that this is another instance where the strength of the $L^2$ orbit equivalence theory allows one to obtain a more general result \cite[Théorème 6.8]{Gab02} in comparison with \cite[Theorem  3.3 (5)]{Luc98b} where a parasitic assumption remains on $Q$ (containing an infinite order element or arbitrarily large finite subgroups) which always holds in a measurable sense. 
In higher degrees the same paradigm is used in \cite[Corollary 1.8]{Sauer-Thom-spectral-seq-L2-2010} (see the proof of Proposition~\ref{prop: L2 Betti poly-free groups}).
With Ab{\'e}rt we proved that $\beta_2^{(2)}(\Out(\FF_n))=0$ for $n\geq 5$ \cite{Abert-Gab}.

The reason for the non-vanishing of $\beta_{2n-3}^{(2)}(\Out(\FF_n))$  and $\beta_{2n-2}^{(2)}(\Aut(\FF_n))$
 boils now down to the existence of subgroups of the form $\FF_2\ltimes \FF_2^{2n-4}$  (resp. $(\FF_2\ltimes \FF_2^{2n-4})\ltimes \FF_n$), to the use of Proposition~\ref{prop: L2 Betti poly-free groups} and to an application of Theorem~\ref{th: top-dim groups} applied to $L=$ the spine of the Culler-Vogtmann space $CV_n$ which is contractible, has dimension $2n-3$, and is equipped with a proper action of $\Out(\FF_n)$
 \cite{CV86} (and its avatar for $\Aut(\FF_n)$).

 Let $(x_1, x_2, \cdots, x_n)$ be a free base of the free group $\FF_n$.
Choose a rank $2$ free subgroup $V\leq\Out(\FF(x_1,x_2))\simeq \GL(2,\Zmath)$ and pick a section $U\leq\Aut(\FF(x_1,x_2))$ of it under 
 $\theta_2$ in the short exact sequence \eqref{eq: Aut to Out}. 
The family of automorphisms  $\Phi(x_1)=\alpha(x_1), \Phi(x_2)=\alpha(x_2), \Phi(x_j)= l_j x_j r_j^{-1}$ ($j\not=1,2$)
for all choices of $(\alpha, l_3,r_3, l_4,r_4, \cdots, l_n,r_n)\in U\times \FF(x_1,x_2)^{2n-4}$ defines a subgroup of $\Aut(\FF_n)$ which is isomorphic to $U\ltimes \FF_2^{2(n-2)}=\FF_2\ltimes \FF_2^{2(n-2)}$  and descends injectively to $\Lambda_n\leq\Out(\FF_n)$ under $\theta_n$. This reproduces an argument from \cite{Bestvina-Kapovich-Kleiner-2002}.
Its pull-back $\widetilde{\Lambda_n}:=\theta_n^{-1}(\Lambda_n)$ is thus isomorphic to $(U\ltimes \FF_2^{2(n-2)})\ltimes \FF_n\simeq\FF_2\ltimes (\FF_2^{2(n-2)}\ltimes \FF_n)$ (the restriction of $\theta_n$ to $\widetilde{\Lambda_n}$ admits a section, thus the splitting).
By Proposition~\ref{prop: L2 Betti poly-free groups}, these poly-free groups satisfy $\beta^{(2)}_{2n-3}(\Lambda_n)=1$ and $\beta^{(2)}_{2n-2}(\widetilde{\Lambda_n})=n-1$. Then apply Theorem~\ref{th: L2 erg-dim groups}.
\CQFD

\begin{proposition}[Poly-free groups] \label{prop: L2 Betti poly-free groups}
\label{cor: beta-d polyfree group}
Consider a group $G=G_n$ obtained by a finite sequence $({G}_i)_{i=1}^{n}$ of extensions 
\begin{equation}1\to {G}_{i}\to {G}_{i+1}\to Q_{i+1}\to 1, \label{eq: extensions of free}\end{equation}
where ${G}_1$ and all the $Q_i$ are finitely generated, non-cyclic free groups.
\\
Then for all $j$ the $\beta_{j}^{(2)}({G}_n)$ vanish except ``in top-dimension'' 
\begin{equation*}
\beta_{n}^{(2)}({G}_n)=\beta_{1}^{(2)}({G}_1)\, \prod_{i=2}^{n} \beta_{1}^{(2)}(Q_i)=(-1)^{n}\chi(G_n).
\end{equation*}
\end{proposition}
Proof: The statement is obtained by induction from the following:
\begin{enumerate}
\item the general results on cohomological/geometric dimension for extensions imply that the geometric dimension of ${G}_n$ is $\leq n$ \cite[Chapter VIII.2]{Bro83};

\item   a result \cite[Theorem  3.3 (5)]{Luc98b}, \cite[Corollary 1.8]{Sauer-Thom-spectral-seq-L2-2010} alluded to above:
{\em Let $1\to N\to \Gamma\to Q\to 1$ be a short exact sequence of infinite groups.
If $\beta_k^{(2)}(N)=0$ for $k=0, 1, \cdots, d-1$ and $\beta_d^{(2)}(N)<\infty$, then $\beta_k^{(2)}(\Gamma)=0$ for $k=0, 1, \cdots d$};

\item the multiplicativity of the Euler characteristic under extensions  \cite[Chapter IX.7]{Bro83} and the coincidence of Euler and $\ell^2$-Euler characteristics \cite[Proposition 0.4]{CG86}:
\begin{eqnarray*}
\chi^{(2)}({G}_{i+1})&=& \chi({G}_{i})\cdot \chi(Q_{i+1})
= (-1)^{i}\beta_{i}^{(2)}({G}_{i}) (-1)  \beta_1^{(2)}(Q_{i+1})
\\
&=& \sum_{j=0}^{\infty} (-1)^{j} \beta_{j}^{(2)}({G}_{i+1})=
(-1)^{i+1}\beta_{i+1}^{(2)}({G}_{i+1}).
\end{eqnarray*}\CQFD
\end{enumerate}

\section{Proof of Theorem~\ref{th: top dim beta 2 of Torelli for Out Fn and Aut Fn} for the Torelli subgroups}
\label{sect: Torelli subgroup}

We continue with the notation of the previous section.
 Pick two elements that generate a free subgroup of rank $2$  in the intersection of 
 the commutator subgroup $[\FF_n,\FF_n]$ with $\FF(x_1,x_2)$, for instance  $u:=[x_1,x_2]$ and $v:=[x_1^{-1},x_2^{-1}]$.
The family of automorphisms  $\Phi(x_1)=x_1, \Phi(x_2)=x_2, \Phi(x_j)= l_j x_j r_j^{-1}$ ($j\not=1,2$)
for all choices of $( l_3,r_3, l_4,r_4, \cdots, l_n,r_n)\in \FF(u,v)^{2n-4}$ defines a subgroup of $\Aut(\FF_n)$ which is isomorphic to $\FF_2^{2n-4}$  and descends injectively under $\theta_n$ (of the exact sequence \eqref{eq: Aut to Out}) to $\Delta_n\leq\mathcal{T}_n\leq\Out(\FF_n)$.
Its pullback $\widetilde{\Delta_n}:=\theta_n^{-1}(\Delta_n)$ is thus a subgroup of $\mathcal{K}_n$ isomorphic to $\FF_2^{2n-4}\ltimes \FF_n$.
Proposition~\ref{prop: L2 Betti poly-free groups} gives $\beta_{2n-4}^{(2)}(\Delta_n)=1$ and $\beta_{2n-3}^{(2)}(\widetilde{\Delta_n})=n-1$.
 The group $\mathcal{T}_n$ has cohomological dimension $2n-4$ \cite{Bestvina-Bux-Margalit-2007}. 
By its general behavior under exact sequences and  $1\to \FF_n\to \mathcal{K}_n\to \mathcal{T}_n\to 1$, the cohomological dimension of $\mathcal{K}_n$ is $ 2n-3$. Then apply Theorem~\ref{th: L2 erg-dim groups}.
\CQFD

\section{Proof of Theorem~\ref{th: top dim L2 pmp equiv rel, geom dim}, measured theoretic version}
\label{sect: proof of R-complex top-dim trick}

Let's consider now the measured theoretic version below of Theorem~\ref{th: top-dim groups}. Theorem~\ref{th: top dim L2 pmp equiv rel, geom dim} will follow directly.
We assume some familiarity with the foundations \cite{Gab02} and refer to this for some background.
\begin{theorem}[Top-dimension $\beta_d^{(2)}$, discrete $\RR$-complex version]
\label{th: top dim L2 R-complex}
Let $(X,\mu)$ be a standard probability measure space and let $\RR$ be a pmp equivalence relation.
Assume $\Sigma$ is $d$-dimensional simplicial discrete $\RR$-complex with vanishing top-dimensional $L^2$-Betti number, $\beta_d^{(2)}(\Sigma,\RR,\mu)=0$.
For any sub-equivalence relation $\SS\leq\RR$ the $L^2$-Betti number of $\Sigma$ seen as a simplicial discrete $\SS$-complex also vanishes in degree $d$, i.e.,  $\beta_d^{(2)}(\Sigma,\SS,\mu)=0$.
\end{theorem}
Proof of Th.~\ref{th: top dim L2 R-complex}:
Recall \cite{Dix69,Gab02} that a measurable bundle $x\mapsto \Sigma_x$ over $(X,\mu)$ of simplicial complexes with uniform bounded geometry delivers an integrated field of $\ell^2$-chain complexes 
$ C_{k}^{(2)}(\Sigma)=\int_X^{\oplus}  C_{k}^{(2)}(\Sigma_{x})\ d\mu(x)$,
and that the field of boundary operators can be integrated into a continuous operator
$\partial_{k}=\int_X^{\oplus}  \left(\partial_{{k},x}:C_{k}^{(2)}(\Sigma_{x})\to C_{{k}-1}^{(2)}(\Sigma_{x})\right)\ d\mu(x)$.

By commutation of the diagram involving the boundary operators and the injective operators induced by inclusion, one gets:
\begin{claim}\label{claim:inject of top dim ker}
Let $\Theta$ and $\Omega$ be measurable  bundles $x\mapsto \Theta_x$ and $x\mapsto \Omega_x$ over $(X,\mu)$ of simplicial complexes both with a bounded geometry. If  
$\Theta\subset \Omega$
then 
\begin{equation}
\ker  \left(\partial_{k} : C_d^{(2)}(\Theta)\to C_{{k}-1}^{(2)}(\Theta)\right)\lhook\joinrel\longrightarrow
\ker  \left(\partial_{k} : C_d^{(2)}(\Omega)\to C_{{k}-1}^{(2)}(\Omega)\right).
\end{equation}
\end{claim}

Recall from \cite[Définition 2.6 and Définition 2.7]{Gab02}) that 
a simplicial discrete (or smooth) $d$-dimensional $\RR$-complex $\Sigma$ 
is an $\RR$-equivariant measurable bundle $x\mapsto \Sigma_x$ of simplicial complexes over $(X,\mu)$
\\
--  that is discrete (the $\RR$-equivariant field of $0$-dimensional cells $\Sigma^{(0)}: x\mapsto \Sigma_x^{(0)}$ admits a Borel fundamental domain); and
\\
-- such that ($\mu$-almost) every fiber $\Sigma_x$ is $\leq d$-dimensional and  $\Sigma_x$ is $d$-dimensional for a non-null set of $x\in X$.

Recall that such an $\RR$-complex is called \defin{uniformly locally bounded (ULB)} if $\Sigma^{(0)}$ admits a finite measure fundamental domain (for its natural fibered measure) and if it admits a uniform bound on the valency of ($\mu$-almost) every vertex  $v\in \Sigma^{(0)}$ (uniform bounded geometry).
Recall the definition of the $L^2$-Betti numbers of the $\RR$-complex $\Sigma$
 \cite[Définition 3.7 and Proposition 3.9]{Gab02}:
\\
Choose any sequence $(\Sigma_{i})_{i}$ of ULB $\RR$-invariant subcomplexes of $\Sigma$ (given by the sequence of bundles $x\mapsto \Sigma_{{i},x}$) which is increasing ($\Sigma_{i}\subset \Sigma_{{i}+1}$) and exhausting ($\cup_{i} \Sigma_{i}=\Sigma$). Let's call such a sequence a \defin{good $\RR$-exhaustion} of $\Sigma$. Let $\mathcal{M}(\RR)$ be the von Neumann algebra of $\RR$.
 The continuous boundary operators ${\partial_{k}}$
 $$0\overset{\partial_0}{\longleftarrow} C_0^{(2)}(\Sigma_{i})\overset{\partial_1}{\longleftarrow} C_1^{(2)}(\Sigma_{i})\overset{\partial_2}{\longleftarrow} 
\cdots \overset{\partial_{k}}{\longleftarrow} C_{k}^{(2)}(\Sigma_{i}) \overset{\partial_{{k}+1}}{\longleftarrow} \cdots,$$
are $\mathcal{M}(\RR)$-equivariant between the Hilbert $\mathcal{M}(\RR)$-modules $C_{k}^{(2)}(\Sigma_{i}) $.

 The reduced $L^2$-homology of $\Sigma_{i}$ is defined as expected as the  Hilbert $\mathcal{M}(\RR)$-module quotient of the kernel by the closure of the image:
\begin{equation*}
\bar{H}_{k}^{(2)}(\Sigma_{i})=\frac{\ker  \left(\partial_{k} : C_{k}^{(2)}(\Sigma_{i})\to C_{{k}-1}^{(2)}(\Sigma_{i})\right)}{\overline{\image}\ \partial_{{k}+1}
\left(\partial_{{k}+1} : C_{{k}+1}^{(2)}(\Sigma_{i})\to C_{k}^{(2)}(\Sigma_{i})\right)}.
\end{equation*}
The inclusions $\Sigma_{i}\subset \Sigma_{j}$ (for ${i}\leq {j}$) induce Hilbert $\mathcal{M}(\RR)$-module operators $C_{k}^{(2)}(\Sigma_{i})\to C_{k}^{(2)}(\Sigma_{j})$ that descend to Hilbert $\mathcal{M}(\RR)$-module operators $\bar{H}_{k}^{(2)}(\Sigma_{i})\overset{J_{{i}, {j}}}{\longrightarrow} \bar{H}_{k}^{(2)}(\Sigma_{j})$. The ${k}$-th $L^2$-Betti number is the double limit of the von Neumann dimension of the closure of the image of these maps:

\begin{equation*}
\beta_{k}^{(2)}(\Sigma,\RR,\mu)=\lim\limits_{{i}\to \infty}\nearrow 
\lim\limits_{
\begin{smallmatrix} {j}\to \infty\\ {j}\geq {i}
\end{smallmatrix}
}\searrow 
\dim_{\mathcal{M}(\RR)}\overline{\image}\ \left(\bar{H}_{k}^{(2)}(\Sigma_{i})\overset{J_{{i}, {j}}}{\longrightarrow} \bar{H}_{k}^{(2)}(\Sigma_{j})\right)   .
\end{equation*}

\medskip
\begin{claim}
\label{claim: beta d non zero iff ker dd non zero}
In the particular case when ${k}=d$ is the top-dimension of $\Sigma$ and $(\Sigma_{i})_{i}$ is a good $\RR$-exhaustion of $\Sigma$, then 
we have the equivalence: $\beta_d^{(2)}(\Sigma,\RR,\mu)>0$ if and only if 
$\ker  \left(\partial_{d} : C_d^{(2)}(\Sigma_{i})\to C_{d-1}^{(2)}(\Sigma_{i})\right)\not=\{0\}$ for a large enough ${i}$.
\end{claim}
Proof: Since $C_{d+1}^{(2)}(\Sigma_{i})=\{0\}$ then 
$\bar{H}_d^{(2)}(\Sigma_{i})=\ker  \left(\partial_{d} : C_d^{(2)}(\Sigma_{i})\to C_{d-1}^{(2)}(\Sigma_{i})\right)
$ for every ${i}$.
Thus by Claim~\ref{claim:inject of top dim ker}
\[\overline{\image}\ \bar{H}_d^{(2)}(\Sigma_{i})\overset{J_{{i}, {j}}}{\longrightarrow} \bar{H}_d^{(2)}(\Sigma_{j})= \ker  \left(\partial_{d} : C_d^{(2)}(\Sigma_{i})\to C_{d-1}^{(2)}(\Sigma_{i})\right).\]
Then
$\beta_{d}^{(2)}(\Sigma,\RR,\mu)=\lim\limits_{{i}\to \infty}\nearrow \dim_{\mathcal{M}(\RR)} \ker  \left(\partial_{d} : C_d^{(2)}(\Sigma_{i})\to C_{d-1}^{(2)}(\Sigma_{i})\right)$.
The claim~\ref{claim: beta d non zero iff ker dd non zero} follows by faithfulness: the property that the von Neumann dimension is non zero if and only if the Hilbert module is non zero.\CQFD

\medskip
The  $d$-dimensional simplicial discrete  $\RR$-complex $\Sigma$ is also an $\SS$-complex with the same properties.
Let $(\Omega_{i})_{i}$ be a good $\RR$-exhaustion of $\Sigma$
 and let $(\Theta_{i})_{i}$ be a similar good $\SS$-exhaustion of $\Sigma$ such that $\Theta_{i}\subset \Omega_{i}$ (one can for instance consider the intersection of a good $\SS$-exhaustion of $\Sigma$ with the good $\RR$-exhaustion $(\Omega_{i})_{i}$).
Assume by contraposition that $\beta_d^{(2)}(\Sigma,\SS,\mu)>0$. It follows from Claim \ref{claim: beta d non zero iff ker dd non zero} that  $\ker  \left(\partial_{d} : C_d^{(2)}(\Theta_{i})\to C_{d-1}^{(2)}(\Theta_{i})\right)\not=\{0\}$ for $\Theta_{i}$ and a large enough ${i}$. Then  the same holds, $\ker  \left(\partial_{d} : C_d^{(2)}(\Omega_{i})\to C_{d-1}^{(2)}(\Omega_{i})\right)\not=\{0\}$, for $\Omega_{i}$ by Claim \ref{claim:inject of top dim ker}. It follows that $\beta_d^{(2)}(\Sigma,\RR,\mu)>0$ by Claim \ref{claim: beta d non zero iff ker dd non zero}.
This completes the proof of Th.~\ref{th: top dim L2 R-complex}.\CQFD

\medskip
As for the proof of Theorem~\ref{th: top dim L2 pmp equiv rel, geom dim}, 
recall from \cite[Définition 3.18]{Gab02} that $\RR$ has \defin{geometric dimension} $\leq d$ if it admits a contractible $d$-dimensional simplicial discrete $\RR$-complex $\Sigma$ (see  \cite[Définition 2.6 and Définition 2.7]{Gab02}).
Recall also the definition of the $L^2$-Betti numbers of $\RR$ \cite[Définition 3.14, Théorème 3.13]{Gab02}: $\beta_k^{(2)}(\RR,\mu):=\beta_k^{(2)}(\Sigma,\RR,\mu)$ where $\Sigma$ is any contractible simplicial discrete $\RR$-complex.
A contractible $d$-dimensional simplicial discrete  $\RR$-complex $\Sigma$ is also an $\SS$-complex with the same properties, so that it can be used to compute the $L^2$-Betti numbers of $\SS$.
Thus Theorem~\ref{th: top dim L2 pmp equiv rel, geom dim} is a specialisation of Theorem~\ref{th: top dim L2 R-complex} when $\Sigma$ is contractible.
\CQFD

\section{Proof of Theorem~\ref{th: L2 erg-dim groups} and Corollary~\ref{cor: ergim and L2}}
\label{sect: proofs with erg dim}

Recall from \cite[Définition 6.4]{Gab02} that a group $\Gamma$ has \defin{ergodic dimension} $\leq d$ if it admits a probability measure preserving \alert{free} action $\Gamma\acting^{\alpha}(X,\mu)$ on some standard space such that the orbit equivalence relation $\RR_{\alpha}$ has geometric dimension $\leq d$.
Equivalently, it admits a $\Gamma$-equivariant  bundle $\Sigma: x\mapsto \Sigma_x$ over $(X,\mu)$ of contractible simplicial complexes of dimension $\leq d$ which is measurable and discrete.
See \cite{Gab02,Gab-erg-dim} for more information on ergodic dimension.

\medskip
Proof of Theorem~\ref{th: L2 erg-dim groups}: Assume $\Gamma$ has ergodic dimension $\leq d$ and that this is witnessed by  $\Gamma\acting^{\alpha}(X,\mu)$ and $\Sigma$, a free pmp $\Gamma$-action and a contractible $d$-dimensional simplicial discrete $\RR_{\alpha}$-complex.
The restriction $\omega$ of the action $\alpha$ to $\Lambda$ being also free, the complex $\Sigma$ computes both the $\ell^2$-Betti numbers of $\Gamma$ and of $\Lambda$; more precisely, 
$\beta_k^{(2)}(\Sigma,\RR_\alpha,\mu)=\beta_k^{(2)}(\Gamma)$ and (considering $\Sigma$ as an $\RR_{\omega}$-complex) $\beta_k^{(2)}(\Sigma,\RR_\omega,\mu)=\beta_k^{(2)}(\Lambda)$ \cite[Corollaire 3.16]{Gab02}.
The $L^2$-Betti numbers of $\Sigma$ vanish strictly above its dimension $d$. 
If moreover $\beta_d^{(2)}(\Sigma,\RR_\alpha,\mu)=\beta_{d}^{(2)}(\Gamma)=0$, then applying Theorem~\ref{th: top dim L2 R-complex} gives $\beta_d^{(2)}(\Sigma,\RR_\omega,\mu)=\beta_d^{(2)}(\Lambda)=0$. \CQFD

\medskip
Proof of Corollary~\ref{cor: ergim and L2}:
By \cite[Theorem 0.2 and Proposition 2.7]{CG86}, all the $\ell^2$-Betti numbers of $\Gamma=\Lambda\times  B$ equal $0$, in particular $\beta_{d+1}^{(2)}(\Gamma)=0$. 
By Theorem~\ref{th: L2 erg-dim groups}, the ergodic dimension of $\Gamma$ is $\geq d+1$.
On the other hand, the ergodic dimension of $B$ is $1$ by Ornstein-Weiss \cite{OW80} and the ergodic dimension of a direct sum is bounded above by the sum of the ergodic dimensions of the factors. \CQFD

\section{Proof of Theorems~\ref{th: d-mfld and top L2 Betti} and~\ref{th: 3-dim manifold and L2 Betti} on manifolds}
\label{sect: proof th manifolds}

Proof of Theorem~\ref{th: d-mfld and top L2 Betti}:
By \cite{Conley-Gaboriau-Marks-Tucker-Drob-2016} the fundamental group $\Gamma=\pi_1(M)$ of a compact connected \alert{aspherical} manifold $M$ of dimension $d\geq 3$ has ergodic dimension $\leq d-1$. Then apply Theorem~\ref{th: L2 erg-dim groups}.\CQFD

\medskip
Any improvement on the ergodic dimension of $\pi_1(M)$ would produce in return a corresponding improvement in Theorem~\ref{th: d-mfld and top L2 Betti}.
\begin{question} \label{quest: erg dim asph mfld groups}
What is the ergodic dimension of the fundamental group of a closed connected \alert{hyperbolic} $d$-manifold $M$? Is it $d/2$ when $d$ is even and $(d+1)/2$ when $d$ is odd?
More generally, is the ergodic dimension of the fundamental group of a closed connected aspherical manifold of dimension $d$ bounded above by $(d+1)/2$?
\end{question}

\medskip
Proof of Theorem~\ref{th: 3-dim manifold and L2 Betti}:
Let $\Gamma$ be  the  fundamental group of a connected compact $3$-dimensional manifold $M$. 
If $M$ is non-orientable, then the fundamental group of its orientation covering $\bar{M}\to M$ has index $2$ in $\pi_1(M)$ so that $\tilde{\Lambda}:=\Lambda\cap \pi_1(\bar{M})$ has index $i=1$ or $i=2$ in $\Lambda$ and $\beta_k^{(2)}(\tilde{\Lambda})=[\Lambda:\tilde{\Lambda}]\, \beta_k^{(2)}(\Lambda)$ for every $k$.
Thus, without loss of generality, one can assume that $M$ is orientable.

Recall that a compact $3$-manifold $M$ is \defin{prime} when every connected sum decomposition $M=N_1\sharp N_2$ is trivial in the sense that either $N_1$ or $N_2\simeq \Smath^2$.
Except for  $\Smath^1\times \Smath^2$, the orientable prime manifolds $M$ are \defin{irreducible}: once the potential boundary spheres have been filled in with $3$-balls (which produces $M'$ and does not change the fundamental group), every embedded $2$-sphere bounds a $3$-ball. 

\begin{theorem}[Kneser-Milnor { \cite{Kneser-1929,Milnor-1962}}]
Let $M^3$ be a connected compact orientable manifold. It can be decomposed as a connected sum (along separating spheres) 
$M=M_1\sharp M_2\sharp \ldots \sharp M_k$  whose pieces $M_j$ are prime; i.e., either are\\
$\left\{
\text{
\begin{tabular}{ll}
--  copies of $\Smath^1\times \Smath^2$ (thus $\pi_1(M_j)\simeq \Zmath$), or \\
-- irreducible manifolds 
$\left\{
\text{
\begin{tabular}{ll}
-- that either have finite $\pi_1$, or   \\
--  $\pi_1(M_j)$ is
 the fundamental group of an \\ aspherical orientable $3$-manifold $M_j'$.
\end{tabular}
}
\right.$
\end{tabular}
}
\right.$
\end{theorem}
It follows that the fundamental group of $M$ decomposes as a free product $\pi_1(M)=\pi_1(M_1)*\pi_1(M_2)*\ldots *\pi_1(M_k)$ of copies of $\Zmath$, of finite groups and of $\pi_1$ of aspherical $3$-manifolds; for these, $\beta_k^{(2)}(\pi_1(M_j)\acting \widetilde{M_j})=\beta_k^{(2)}(\pi_1(M_j))$.

The second $\ell^2$-Betti number of the fundamental group of a compact connected orientable irreducible non-exceptional aspherical $3$-manifold vanishes \cite[Theorem~0.1]{LL95}: 
$\beta_2^{(2)}(\pi_1(M_j)\acting \widetilde{M_j})=0$.
By the work of Perelman and his proof of Thurston's geometrisation conjecture \cite{Perelman-2002-entropy-Ricci,Perelman-2003-Ricci-surgery} (see also \cite{Kleiner-Lott-2008-Notes-Perelman,BBBMP-2010}) exceptional manifolds do not exist.

It follows (by the $\ell^2$-version of Mayer-Vietoris \cite{CG86}) that  $\beta_2^{(2)}(\pi_1(M))=0$ for every connected compact $3$-manifold $M$.

The above free product decomposition implies that $\pi_1(M)$ has virtual geometric dimension $\leq 3$.
Moreover by  \cite{Conley-Gaboriau-Marks-Tucker-Drob-2016}, $\pi_1(M)$ has ergodic dimension $\leq 2$. 
Theorem~\ref{th: 3-dim manifold and L2 Betti} then follows from Theorem~\ref{th: L2 erg-dim groups}.
When $\Lambda$ is infinite, $\beta_0^{(2)}(\Lambda)=0$ and $\chi^{(2)}(\Lambda)=\sum_{k}(-1)^k\beta_k^{(2)}(\Lambda)=-\beta_1^{(2)}(\Lambda)\in [-\infty, 0]$.

We now give an alternative argument avoiding the use of the unpublished article \cite{Conley-Gaboriau-Marks-Tucker-Drob-2016}.
If $M$ is an aspherical orientable $3$-manifold with boundary, then its fundamental group has geometric dimension $\leq 2$.
Otherwise, by Thurston's geometrization conjecture (now established), an aspherical orientable $3$-manifold can be decomposed along a disjoint union of embedded tori into pieces which carry a geometric structure. This delivers a further decomposition of its fundamental group as a graph of groups with edge groups isomorphic to $\Zmath^2$.
The fundamental group $\pi_1(M)$ eventually follows decomposed as a graph of groups with edge groups isomorphic to either $\{1\}$ or $\Zmath^2$. The vertex groups $\Gamma_i$ have ergodic dimension $\leq 2$. More precisely, the $\Gamma_i$ are either
\\
--  amenable: they have ergodic dimension $\leq 1$ by \cite{OW80}; or
\\
-- 
a cocompact lattice in the isometry group of one of the Thurston's geometries:  
when it is non-amenable, $\Gamma_i$ is measure equivalent with some non-cocompact  lattice $\Gamma_i'$ in the isometry group of  $\Hmath^3, \Hmath^2\times \Rmath$ or $\widetilde{\PSL(2,\Rmath)}$ ($\Gamma_i'$ has geometric dimension $\leq 2$). Then $\Gamma_i$ has ergodic dimension $\leq 2$ \cite[Proposition 6.5]{Gab02};~or
\\
-- the fundamental group of an aspherical complex of dimension $\leq 2$  (by a deformation retraction of a $3$-dimensional manifold with boundary).

By Mayer-Vietoris \cite{CG86} and by triviality of  $\beta_p^{(2)}$ ($p= 1, 2$)  for amenable groups, $\beta_2^{(2)}(\pi_1(M))$ equals to the sum of the $\beta_2^{(2)}(\Gamma_i)$ of the vertex groups. 
Since $\beta_2^{(2)}(\pi_1(M))=0$, all the vertex groups $\Gamma_i$ satisfy $\beta_2^{(2)}(\Gamma_i)=0$. 
The same holds for their subgroups by Theorem~\ref{th: L2 erg-dim groups}.

A subgroup $\Lambda$ of $\pi_1(M)$ decomposes, by Bass-Serre theory \cite{Ser77}, as a graph of groups whose edge groups are subgroups of $\Zmath^2$ and vertex groups are subgroups of the $\Gamma_i$. Again by  Mayer-Vietoris, $\beta_2^{(2)}(\Lambda)=0$.
\CQFD

\newcommand{\etalchar}[1]{$^{#1}$}
\def\cprime{$'$}


\begin{thebibliography}{BBB{\etalchar{+}}10}

\bibitem[AG20]{Abert-Gab}
M.~Ab{\'e}rt and D.~Gaboriau.
\newblock Higher dimensional cost and profinite actions.
\newblock in preparation, 2022.

\bibitem[Ati76]{Ati76}
{}M. Atiyah.
\newblock Elliptic operators, discrete groups and von {N}eumann algebras.
\newblock In {\em Colloque ``Analyse et Topologie'' en l'Honneur de Henri
  Cartan (Orsay, 1974)}, pages 43--72. Ast\'erisque, SMF, No. 32--33. Soc.
  Math. France, Paris, 1976.

\bibitem[Bar16]{Bartholdi-2016-homol-out-F7}
L. Bartholdi.
\newblock The rational homology of the outer automorphism group of {$\FF_7$}.
\newblock {\em New York J. Math.}, 22:191--197, 2016.

\bibitem[Bau63]{Baumslag-1963-Autom-gps-resid-finite}
G. Baumslag.
\newblock Automorphism groups of residually finite groups.
\newblock {\em J. London Math. Soc.}, 38:117--118, 1963.

\bibitem[BBB{\etalchar{+}}10]{BBBMP-2010}
L.Bessi\`eres, G. Besson, M. Boileau, S. Maillot, and
  J. Porti.
\newblock {\em Geometrisation of 3-manifolds}, volume~13 of {\em EMS Tracts in
  Mathematics}.
\newblock European Mathematical Society (EMS), Z\"{u}rich, 2010.

\bibitem[BBM07]{Bestvina-Bux-Margalit-2007}
M. Bestvina, K.-U. Bux, and D. Margalit.
\newblock Dimension of the {T}orelli group for {${\rm Out}(F_n)$}.
\newblock {\em Invent. Math.}, 170(1):1--32, 2007.

\bibitem[BG20]{Bar-Gab}
L.~Bartholdi and D.~Gaboriau.
\newblock Around the homology of ${\Out(\FF_n)}$.
\newblock in preparation, 2022.

\bibitem[BKK02]{Bestvina-Kapovich-Kleiner-2002}
M. Bestvina, M. Kapovich, and B. Kleiner.
\newblock Van {K}ampen's embedding obstruction for discrete groups.
\newblock {\em Invent. Math.}, 150(2):219--235, 2002.

\bibitem[Bor85]{Bor85}
A.~Borel.
\newblock The ${L}\sp 2$-cohomology of negatively curved {R}iemannian symmetric
  spaces.
\newblock {\em Ann. Acad. Sci. Fenn. Ser. A I Math.}, 10:95--105, 1985.

\bibitem[Bro82]{Bro83}
{}K. Brown.
\newblock {\em Cohomology of groups}.
\newblock Springer-Verlag, New York, 1982.

\bibitem[BV06]{Bridson-Vogtmann-2006}
M.~R. Bridson and K. Vogtmann.
\newblock Automorphism groups of free groups, surface groups and free abelian
  groups.
\newblock In {\em Problems on mapping class groups and related topics},
  volume~74 of {\em Proc. Sympos. Pure Math.}, pages 301--316. Amer. Math.
  Soc., Providence, RI, 2006.

\bibitem[BV19]{Borinsky-Vogtmann-2019}
M. {Borinsky} and K. {Vogtmann}.
\newblock {The Euler characteristic of $\operatorname{Out}(F_n)$}.
\newblock {\em {Comment. Math. Helv.}}, 95(4):703--748, 2020.

\bibitem[CG86]{CG86}
{}J. Cheeger and {}M. Gromov.
\newblock ${L}\sb 2$-cohomology and group cohomology.
\newblock {\em Topology}, 25(2):189--215, 1986.

\bibitem[CGMT]{Conley-Gaboriau-Marks-Tucker-Drob-2016}
C.~T. Conley, {}D. Gaboriau, A.~S. Marks, and R.~D.
  {Tucker-Drob}.
\newblock One-ended spanning subforests and treeability of groups.
\newblock {\em preprint}.

\bibitem[Con79]{Connes-79}
A.~Connes.
\newblock Sur la th\'eorie non commutative de l'int\'egration.
\newblock In {\em Alg\`ebres d'op\'erateurs (S\'em., Les Plans-sur-Bex, 1978)},
  pages 19--143. Springer, Berlin, 1979.

\bibitem[CV86]{CV86}
M.~Culler and K.~Vogtmann.
\newblock Moduli of graphs and automorphisms of free groups.
\newblock {\em Invent. Math.}, 84(1):91--119, 1986.

\bibitem[Dix69]{Dix69}
{}J. Dixmier.
\newblock {\em Les alg\`ebres d'op\'erateurs dans l'espace hilbertien
  (alg\`ebres de von {N}eumann)}.
\newblock Gauthier-Villars \'Editeur, Paris, 1969.
\newblock Deuxi\`eme \'edition, revue et augment\'ee, Cahiers Scientifiques,
  Fasc. XXV.

\bibitem[Dod79]{Dodziuk-1979-harmon-sym-mflds}
J. Dodziuk.
\newblock {$L^{2}$} harmonic forms on rotationally symmetric {R}iemannian
  manifolds.
\newblock {\em Proc. Amer. Math. Soc.}, 77(3):395--400, 1979.

\bibitem[Eck00]{Eck00}
{}B. Eckmann.
\newblock Introduction to $l\sb 2$-methods in topology: reduced $l\sb
  2$-homology, harmonic chains, $l\sb 2$-{B}etti numbers.
\newblock {\em Israel J. Math.}, 117:183--219, 2000.
\newblock Notes prepared by Guido Mislin.

\bibitem[EG57]{Eilenberg-Ganea-1957}
S. Eilenberg and T. Ganea.
\newblock On the {L}usternik-{S}chnirelmann category of abstract groups.
\newblock {\em Ann. of Math. (2)}, 65:517--518, 1957.

\bibitem[Gab02]{Gab02}
D.~Gaboriau.
\newblock Invariants {$L\sp 2$} de relations d'\'equivalence et de groupes.
\newblock {\em Publ. Math. Inst. Hautes \'Etudes Sci.}, 95:93--150, 2002.

\bibitem[Gab20]{Gab-erg-dim}
D.~Gaboriau.
\newblock On the ergodic dimension.
\newblock in preparation, 2022.

\bibitem[GL09]{Gaboriau-Lyons-2009}
D.~Gaboriau and R.~Lyons.
\newblock A measurable-group-theoretic solution to von {N}eumann's problem.
\newblock {\em Invent. Math.}, 177(3):533--540, 2009.

\bibitem[Gro75]{Grossman-74}
E.~K. Grossman.
\newblock On the residual finiteness of certain mapping class groups.
\newblock {\em J. London Math. Soc. (2)}, 9:160--164, 1974/75.

\bibitem[Hat95]{Hatcher-1995-Homol-stab-aut-free-gps}
A. Hatcher.
\newblock Homological stability for automorphism groups of free groups.
\newblock {\em Comment. Math. Helv.}, 70(1):39--62, 1995.

\bibitem[Kam19]{Kammeyer-book-intro-L2}
{}H. Kammeyer.
\newblock {\em Introduction to {$\ell^2$}-invariants}.
\newblock Springer-Verlag, 2019.

\bibitem[Kid08]{Kida-08}
Y.~Kida.
\newblock The mapping class group from the viewpoint of measure equivalence
  theory.
\newblock {\em Mem. Amer. Math. Soc.}, 196(916):viii+190, 2008.

\bibitem[KL08]{Kleiner-Lott-2008-Notes-Perelman}
B. Kleiner and J. Lott.
\newblock Notes on {P}erelman's papers.
\newblock {\em Geom. Topol.}, 12(5):2587--2855, 2008.

\bibitem[Kne29]{Kneser-1929}
H. Kneser.
\newblock Geschlossene Fl{\"a}chen in dreidimensionalen Mannigfaltigkeiten.
\newblock {\em Jahresbericht der Deutschen Mathematiker-Vereinigung},
  38:248--259, 1929.

\bibitem[LL95]{LL95}
{}J. Lott and {}W. L{\"u}ck.
\newblock ${L}\sp 2$-topological invariants of $3$-manifolds.
\newblock {\em Invent. Math.}, 120(1):15--60, 1995.

\bibitem[L{\"u}c94]{Luck-1994-approx}
W.~L{\"u}ck.
\newblock Approximating {$L^2$}-invariants by their finite-dimensional
  analogues.
\newblock {\em Geom. Funct. Anal.}, 4(4):455--481, 1994.

\bibitem[L{\"u}c98]{Luc98b}
{}W. L{\"u}ck.
\newblock Dimension theory of arbitrary modules over finite von {N}eumann
  algebras and ${L}\sp 2$-{B}etti numbers. {I}{I}. {A}pplications to
  {G}rothendieck groups, ${L}\sp 2$-{E}uler characteristics and {B}urnside
  groups.
\newblock {\em J. Reine Angew. Math.}, 496:213--236, 1998.

\bibitem[L{\"u}c02]{Luc02}
W.~L{\"u}ck.
\newblock {\em {$L\sp 2$}-invariants: theory and applications to geometry and
  {$K$}-theory}, volume~44.
\newblock Springer-Verlag, Berlin, 2002.

\bibitem[Mag35]{Magnus-1935}
W. Magnus.
\newblock \"{U}ber {$n$}-dimensionale {G}ittertransformationen.
\newblock {\em Acta Math.}, 64(1):353--367, 1935.

\bibitem[Mil62]{Milnor-1962}
J.~Milnor.
\newblock A unique decomposition theorem for {$3$}-manifolds.
\newblock {\em Amer. J. Math.}, 84:1--7, 1962.

\bibitem[Nie24]{Nielsen-1924}
J. Nielsen.
\newblock Die {I}somorphismengruppe der freien {G}ruppen.
\newblock {\em Math. Ann.}, 91(3-4):169--209, 1924.

\bibitem[Oha08]{Ohashi-2008-rat-homol-Out-Fn}
R. Ohashi.
\newblock The rational homology group of {${\rm Out}(F_n)$} for {$n\leq 6$}.
\newblock {\em Experiment. Math.}, 17(2):167--179, 2008.

\bibitem[OW80]{OW80}
{}D. Ornstein and {}B. Weiss.
\newblock Ergodic theory of amenable group actions. {I}. {T}he {R}ohlin lemma.
\newblock {\em Bull. Amer. Math. Soc. (N.S.)}, 2(1):161--164, 1980.

\bibitem[{Per}02]{Perelman-2002-entropy-Ricci}
G. {Perelman}.
\newblock {The entropy formula for the Ricci flow and its geometric
  applications}.
\newblock {\em arXiv Mathematics e-prints}, page math/0211159, Nov 2002.

\bibitem[{Per}03]{Perelman-2003-Ricci-surgery}
G. {Perelman}.
\newblock {Ricci flow with surgery on three-manifolds}.
\newblock {\em arXiv Mathematics e-prints}, page math/0303109, Mar 2003.

\bibitem[Ser77]{Ser77}
J.-P. Serre.
\newblock {\em Arbres, amalgames, {SL}$\sb{2}$}.
\newblock Ast\'erisque, No. 46. Soci\'et\'e Math\'ematique de France, Paris,
  1977.
  
\bibitem[ST10]{Sauer-Thom-spectral-seq-L2-2010}
R.~Sauer and A.~Thom.
\newblock A spectral sequence to compute {$L^2$}-Betti numbers of groups and
  groupoids.
\newblock {\em J. Lond. Math. Soc. (2)}, 81(3):747--773, 2010.

\bibitem[Vog06]{Vogtmann-2006-ICM-survey}
K. Vogtmann.
\newblock The cohomology of automorphism groups of free groups.
\newblock In {\em International {C}ongress of {M}athematicians. {V}ol. {II}},
  pages 1101--1117. Eur. Math. Soc., Z\"{u}rich, 2006.

\end{thebibliography}
\end{document}